%% file: NewDraft.tex
\documentclass[a4paper,UKenglish,cleveref, autoref, thm-restate,numberwithinsect]{lipics-v2021}

\hideLIPIcs  

\title{Order-Invariant CFO Contains Local-Order-Invariant FO} 

\author{Derek {R. Aoki}}{Proof School, San Francisco, United States \and \url{http://www.derekaoki.com} }{derekaoki@acm.org}{https://orcid.org/0000-0003-0352-1963}{}

\authorrunning{D. R. Aoki} 

\Copyright{Derek R. Aoki} 



\category{} 

\relatedversion{} 

\acknowledgements{I would like to thank Noah Schweber for helpful conversations during the writing of this article.}

\nolinenumbers 

\input{preamble}

\begin{document}

\maketitle

\begin{abstract}
    Order-invariant cluster first-order logic (CFO) and local-order-invariant first-order logic are two fragments of order-invariant first-order logic, both motivated by the study of invariant logics on classes of bounded degree. Our main result is that order-invariant CFO contains local-order-invariant FO on arbitrary classes of finite structures, settling the conjecture of Ghasemi and Grange that order-invariant CFO contains plain FO. By combining this with a collapse result of Ghasemi and Grange, we obtain the equivalence of FO, local-order-invariant FO, and order-invariant CFO on classes of finite colored graphs of bounded degree. The proof makes use of Gaifman normal form, which allows our translation from local-order-invariant FO to order-invariant CFO to be effective. Combined with the fixed-parameter tractability of model-checking for order-invariant CFO on bounded-degree colored graph classes, this shows that model-checking for local-order-invariant FO on bounded-degree colored graph classes is fixed-parameter tractable. To generalize our collapse and model-checking results beyond colored graphs, we use a method based on first-order interpretations which also yields a collapse of order-invariant FO to FO when restricted to graph classes of bounded shrub-depth. Finally, we give an application to the database-motivated frozen-witness logic of Abiteboul, Simon, and Vianu by showing that local-order-invariant FO contains deterministic frozen-witness FO, and thus that deterministic frozen-witness FO collapses to FO on classes of finite structures of bounded degree.
\end{abstract}

\section{Introduction}

Finite model theory is the study of the queries on databases that can be expressed in formal languages. In particular, one tends to be interested in which queries are expressible in a given logic, and how difficult it is to compute whether a database satisfies a query defined in a given logic. In standard logics like first-order logic, fixed-point logics, or second-order logic, one can only refer to the relations of the structure itself. However, it is natural to want to make use of information like the order the data in the database is stored in, so long as the query evaluates in the same way on a database regardless of the particular order chosen. We call such a query, that is, one that uses an order but does not depend on the particular choice of order, order-invariant.

In the case of structures of arbitrary (including infinite) size, it is known that order-invariant first-order logic is equivalent in expressiveness to plain $\FO$, via the Craig Interpolation Theorem \cite{craig1957}. However, there exist properties of finite structures definable in order-invariant $\FO$ which are not definable in $\FO$, the first of which was discovered by Gurevich  and published in \cite{ebbinghausFlum2005}. In fact, it has been shown that the weaker fragment successor-invariant $\FO$, in which formulae can only access the successor relation induced by the order in which the data is stored, can express properties of finite structures that cannot be expressed in $\FO$ \cite{rossman2007}. Further non-collapse results include those for epsilon-invariant $\FO$ and deterministic frozen-witness $\FO$ \cite{otto2000}, as well as for local-order-invariant $\FO$ \cite{aoki2025}.

Based on this, it is natural to ask whether invariant logics collapse on finite structures or on particular classes of finite structures. For instance, it has been proved that successor-invariant $\FO$ and $\FO$ express the same queries on classes of bounded degree \cite{grange2021}, and that order-invariant $\FO$ collapses on trees, strings, and hollow trees \cite{benediktSegoufin2009,grangeSegoufin2020}. The question of whether full order-invariant $\FO$ collapses on classes of bounded degree is significant and open. However, a substantial fragment of order-invariant $\FO$, known as order-invariant cluster first-order logic ($\CFO$), has recently been introduced and shown to be contained in $\FO$ on classes of colored graphs of bounded degree \cite{ghasemigrange}. A second question concerns the complexity of model-checking for these invariant logics, and it has been shown for instance that model-checking for order-invariant $\CFO$ is FPT \cite{ghasemigrange}. 

An invariant logic of particular interest to us is local-order-invariant logic. A local order on a graph is a ternary relation $\prec(x,y,z)$ such that when $x$ is fixed, $\prec(x,y,z)$ linearly orders the neighbors of $x$. On structures over arbitrary vocabularies, a local order is a local order on the Gaifman graph of the structure. Such local orders have been considered in the context of finite model theory prior to being considered in the context of invariant logics \cite{benediktSegoufin2009,courcelle1996}. Local-order-invariant $\FO$ was first considered in \cite{lindellTowsnerWeinstein2025}, and a number of questions related to this work were asked by Weinstein in \cite{Weinstein2018InvariantDefinability}. Generalizations of these questions were answered in \cite{aoki2025}.

\textbf{Our Contributions}

\begin{enumerate}
    \item Our main result is that every property expressible in local-order-invariant $\FO$ is expressible in order-invariant $\CFO$. Our translation uses Gaifman normal form, as introduced in \cite{gaifman1982}, which makes it effective. Thus, our translation works on arbitrary classes of finite structures. This affirms in particular the conjecture of Ghasemi and Grange in \cite{ghasemigrange} that order-invariant $\CFO$ contains $\FO$. An important corollary is that $\FO$ and local-order-invariant $\FO$ have the same expressiveness on classes of bounded degree. In \cite{aoki2025}, this result is proved in particular for the full class of structures over a vocabulary $\sigma$ of degree at most $d$, and our result generalizes this to arbitrary classes of bounded degree. Another corollary is that, by the effectiveness of the translation to order-invariant $\CFO$, the model-checking problem for local-order-invariant $\FO$ is FPT. All of these results are proved first for colored graphs. 
    \item To generalize our main result, its corollaries, and the collapse result of \cite{ghasemigrange} from classes of colored graphs of bounded degree to classes of bounded degree over arbitrary vocabularies, we use a certain class of first-order interpretations. We believe that this method of investigating invariant logics via interpretations has been under-explored, as it allows us, among several other things, to straightforwardly combine several results from \cite{chenFlum2020,eickmeyerElberfeldHarwath2017,ohlmannPilipczukPrzybyszewskiTorunczyk2023} to obtain the collapse of order-invariant $\FO$ to $\FO$ when restricted to graph classes of bounded shrub-depth.
    \item It is proved in \cite{aoki2025} that epsilon-invariant $\FO$ is contained in local-order-invariant $\FO$ on classes of bounded degree, but the proof is flawed in that it accounts for neither the extensionality nor the nesting of epsilon terms. However, the basic idea of the proof is usable. Deterministic frozen-witness $\FO$, introduced and studied in \cite{abiteboulHullVianu1995,abiteboulSimonVianu1990,AbiteboulVianu1991Nondeterminism}, is a weaker logic than epsilon-invariant $\FO$. It is of interest in the study of non-determinism in databases. Abiteboul, Simon, and Vianu originally conjectured in \cite{abiteboulSimonVianu1990} that it should collapse to $\FO$, and this conjecture was refuted in \cite{otto2000}. We use Gaifman normal form again to show the containment of deterministic frozen-witness $\FO$ in local-order-invariant $\FO$ on finite structures, the collapse of deterministic frozen-witness $\FO$ on classes of finite structures of bounded degree, and the fixed-parameter tractability of the model-checking problem for deterministic frozen-witness $\FO$ on classes of finite structures of bounded degree. This recovers the collapse conjectured by Abiteboul, Simon, and Vianu on classes of bounded degree.
\end{enumerate}

In \Cref{sec:preliminaries}, we give the necessary background. In \Cref{sec:maincontainment}, we prove our titular containment theorem and its corollaries in the setting of colored graphs. In \Cref{sec:interpretations}, we extend the main containment and related results to arbitrary vocabularies, and prove the collapse for bounded shrub-depth. In \Cref{sec:frozenwitness}, we prove the deterministic frozen-witness containment and its corollaries.

\section{Preliminaries}\label{sec:preliminaries}

From here, $\sigma$ will be a finite relational vocabulary. We denote by $\mfF_\sigma$ the class of finite structures over $\sigma$, and $\mfK_\sigma$ the class of arbitrary structures over $\sigma$.

\subsection{Locality}

For any $X\subseteq A$, we write $A\harpoon X$ for the restriction of $A$ to $X$; i.e., the substructure of $A$ with domain $X$ such that for any tuple $\vec a$ of elements of $X$ and any relation $R\in \sigma$ of arity $|\vec a|$, we have $R(\vec a)$ in $A\harpoon X$ if and only if $R(\vec a)$.

\begin{definition}
    Let $A$ be in $\mfK_\sigma$. We define the Gaifman graph of $A$, denoted by $\mcG(A)$, such that the vertices of $\mcG(A)$ are the elements of $A$, and such that there are edges between two vertices $a,a'\in A$ if and only if there exist a tuple $\vec a$ of elements of $A$ and a relation $R\in \sigma$ of arity $|\vec a|$ such that $a,a'\in \vec a$ and $R(\vec a)$ in $A$. We now recall several topological graph notions as they relate to Gaifman graphs.
    \begin{enumerate}
        \item For $a,a'\in A$, we define the distance $\delta^A(a,a')$ as the least number of edges in a path in $\mcG(A)$ from $a$ to $a'$, if such a path exists, and $\infty$ otherwise. 
        \item For $a\in A$, we denote by $B_r^A(a)$ the ball of radius $r$ around $a$, $\{a'\in A:\delta^A(a,a')\leq r\}$ . We denote by $S_r^A(a)$ the sphere of radius $r$ around $a$, $B_r^A(a)\setminus B_{r-1}^A(a)$.
        \item For $a\in A$, we denote by $N_r^A(a)$ the structure over $\sigma\cup\{c\}$ which is identical to the structure $A\harpoon B_r^A(a)$ except that it has the additional constant symbol $c$ which is equal to $a$. We call this the neighborhood of radius $r$ around $a$. We say that the $r$-neighborhood type of $a\in A$, denoted $\tp_r^A(a)$, is the class $\{B\in \mfK_\sigma:B\cong N_r^A(a)\}$.
    \end{enumerate}
    We can generalize the above definitions of balls, spheres, and neighborhoods to arbitrary tuples $\vec a$ of elements of $A$ by taking $B_r^A(\vec a)=\bigcup_{a\in \vec a}B_r^A(a)$.
\end{definition}

\begin{definition}
    We define $\mfF_\sigma^d$ as the structures of degree at most $d$ over $\sigma$, i.e. the elements $A$ of $\mfF_\sigma$ such that every vertex in $\mcG(A)$ has degree at most $d$. Similarly, we define $\mfK_\sigma^{<\omega}$ as the locally finite structures over $\sigma$, i.e. elements $A$ of $\mfK_\sigma$ such that every vertex in $\mcG(A)$ has finitely many neighbors.
\end{definition}

\subsection{Invariant Logics}

\begin{definition}
    Let $\varphi$ be a sentence in the vocabulary $\sigma\cup \{<\}$. We say that $\varphi$ is order-invariant on $\mfK\subseteq \mfK_\sigma$ if, for any $A\in \mfK$ and any linear orders $<_1,<_2$ on $A$, we have $(A,<_1)\models \varphi$ if and only if $(A,<_2)\models \varphi$. If $\varphi$ is order-invariant, we say that $A\models \varphi$ if $(A,<)\models \varphi$ for any linear order $<$.
\end{definition}

\begin{definition}
    We say that a ternary relation $\prec(x,y,z)$ on a graph $A$ is a local order if, when $x$ is fixed, $\prec(x,y,z)$ is a linear order on $S_1^A(x)$. This definition generalizes to general structures by instead considering the Gaifman graph $\mcG(A)$. We define local-order-invariance of a sentence $\varphi$ in $\sigma\cup\{\prec\}$ in the natural way.
\end{definition}

The following is the main result of \cite{aoki2025}. A generalization is claimed for arbitrary classes of bounded degree and used to prove the claimed collapse of epsilon-invariant $\FO$ to $\FO$ on classes of bounded degree, but this generalization does not follow from the proof techniques of the paper because the result of \cite{lindellTowsnerWeinstein2025} used does not apply in general (in particular, it has the additional hypothesis that the query must be invariant on an ambient substructure-closed class). 

\begin{theorem}[Aoki]\label{thm:aokicollapse}
    Let $d\in\omega$, and let $\varphi$ be a local-order-invariant sentence over $\mfF_\sigma^d$. Then there is a first-order sentence $\psi$ such that for any $A\in \mfF_\sigma^d$, $A\models \psi$ if and only if $A\models \varphi$.
\end{theorem}

\subsection{Cluster First-Order Logic}

Before describing cluster first-order logic on ordered graphs rigorously, we give a brief informal description. Note that, for the following discussion and definition, the graphs in question will be simple (i.e., undirected and without self-edges). First, we fix an infinite alphabet $\Sigma$. The variables of CFO are indexed as $x_w^i$, with $i\in \omega$ and $w\in \Sigma^*$, and these variables are grouped into clusters based on whether they share the index $w$. Introducing a variable $x_w^j$ with $j>0$ via a quantifier is guarded by adjacency with the previous elements of the cluster, in the sense that there must exist some $i<j$ such that $E(x_w^i,x_w^j)$. Also, the relations $E$, $=$, and $<$ can only be used to compare variables of the same cluster, with the exception that $=$ and $<$ may be used to compare $x_w^i$ with $x_{w\alpha}^0$ for any $i$. Now for the rigorous definition.

\begin{definition}
    Fix an infinite alphabet $\Sigma$, and let $S\subseteq \Sigma^*\times \omega$ be a finite set of indices. We say that $S$ is valid if it is downward closed, in the sense that for any $w\in \Sigma^*$, $\alpha\in \Sigma$, and $i\in \omega$,
    \begin{enumerate}
        \item for every $j<i$, $(w,i)\in S$ implies $(w,j)\in S$, and
        \item $(w\alpha,0)\in S$ implies $(w,0)\in S$. 
    \end{enumerate}
    If $S$ is a valid set of indices, we write $X_S$ for the set $\{x_w^i:(w,i)\in S\}$. Letting $\sigma$ be a vocabulary $\{E,P_1,\dots,P_n\}$, where $E$ is binary and each $P_i$ is unary, we define $\CFO_S[k]$ by induction on $k$.
    \begin{enumerate}
        \item For every $(w,i)\in S$, $j<i$, $l\leq n$, $\alpha\in \Sigma$ such that $(w\alpha,0)\in S$, and $\sim$ in $\{=,<,>\}$, then $\CFO_S[0]$ contains the atomic formulae $P_l(x_w^i)$, $E(x_w^i,x_w^j)$, $x_w^i\sim x_w^j$, and $x_w^i\sim x_{w\alpha}^0$. Also, $\CFO_S[0]$ is closed under Boolean combinations.
        \item $\CFO_S[k+1]$ contains the formulae
        \begin{enumerate}
            \item $\exists x_\varepsilon^0 \varphi$, where $\varphi\in \CFO_{\{(\varepsilon,0)\}}$ (assuming that $S=\emptyset$);
            \item $\exists x_{w\alpha}^0 \varphi$, where $(w,0)\in S$, $\alpha\in \Sigma$, $(w\alpha,0)\not\in S$, and $\varphi\in \CFO_{S\cup \{(w\alpha,0)\}}[k]$;
            \item $\exists x_w^i E(x_w^i,x_w^j)\land \varphi$, where $(w,j)\in S$, $i=\min\{l:(w,l)\not\in S\}$, and $\varphi\in \CFO_{S\cup \{(w,i)\}}[k]$. 
        \end{enumerate}
        Also, $\CFO_S[k+1]$ is closed under Boolean combinations.
    \end{enumerate}
    We write $\CFO_S=\bigcup_{k\in \omega}\CFO_S[k]$, and we write $\CFO$ for $\CFO_\emptyset$. Because $\CFO_S$ is a fragment of first-order logic, its semantics are very natural. That is, for $\varphi\in \CFO_S$, $A\in \mfF_\sigma$, a linear order $<$ on $A$, and a valuation $\nu:X_S\to A$, we have $(A,<,\nu)\models \varphi$ (or $(A,<)\models \varphi$ if $S=\emptyset$) if this holds when $\varphi$ is treated as a first-order formula with free variables in $X_S$. Order-invariance of a $\CFO$ sentence is defined in the obvious way.
\end{definition}

Though $\CFO$ was defined originally in \cite{ghasemigrange} only for (ordered) colored graphs, there is a natural extension to (ordered) structures over any vocabulary. To do this, we first note the following straightforward lemma.

\begin{lemma}\label{lem:gaifmanedgedefinability}
    For every vocabulary $\sigma$, there exists a formula $E_G(x,y)$ such that for any $a,a'$, we have $E_G(a,a')$ in $A$ if and only if $a$ and $a'$ are adjacent in $\mcG(A)$.
\end{lemma}

\begin{definition}
    Define $\CFO$ over an arbitrary vocabulary $\sigma$ identically to $\CFO$ over colored graphs, except
    \begin{enumerate}
        \item $\CFO_S[0]$ contains $R(\vec x)$ where the elements of $\vec x$ are all of the form $x_w^i$ for some $w$ with $(w,i)\in S$, and
        \item we guard cluster quantification by $\exists x_w^i E_G(x_w^i,x_w^j)\land \varphi$, where $E_G$ is defined as above.
    \end{enumerate}
\end{definition}

\begin{theorem}[Ghasemi--Grange]\label{thm:CFOcollapseoncoloredgraphs}
    Let $\sigma$ be a colored graph vocabulary as above, and let $\mfK$ be a class of finite colored graphs over $\sigma$ of degree at most $d$. Let $\varphi$ be an order-invariant $\CFO$ sentence on $\mfK$. Then $\varphi$ is equivalent over $\mfK$ to some first-order sentence $\psi$ in the vocabulary $\sigma$. In other words, order-invariant $\CFO$ is contained in $\FO$ on colored graph classes of bounded degree.
\end{theorem}

The following was explicitly conjectured in \cite{ghasemigrange}.

\begin{conjecture}[Ghasemi--Grange]\label{conj:CFOcontainsFO}
    Let $\mfK$ be a class of finite colored graphs over $\sigma$, and let $\varphi$ be a first-order sentence. Then $\varphi$ is equivalent on $\mfK$ to some order-invariant $\CFO$ sentence $\psi$. In other words, $\FO$ is contained in order-invariant $\CFO$ on arbitrary classes of finite colored graphs.
\end{conjecture}

\section{Main Containment Theorem}\label{sec:maincontainment}

In this section, we will only work with colored simple graphs over a vocabulary $\sigma$.

We will now introduce Gaifman normal form and several relevant definitions. While it is true that in this subsection we are only working with colored simple graphs, the following applies to an arbitrary relational vocabulary.

For $r\in \omega$ and a first-order formula $\varphi(\vec x)$, we write $\varphi^{(r)}$ for the formula which is identical to $\varphi$ except the quantifiers are bounded to $B_r^A(\vec x)$. We call $\varphi^{(r)}$ an $r$-local formula. As a consequence of \Cref{lem:gaifmanedgedefinability}, there exist formulae which express $\delta(x,y)\leq r$ and $\delta(x,y)>r$ for every $r$. So we can instead phrase $\varphi^{(r)}$ as the formula that results when $\exists y \psi$ is replaced by $\exists y(\bigvee_i\delta(x_i,y)\leq r\land \psi)$, and similarly for universal quantifiers. We define basic local sentences as sentences of the form
\begin{align*}
    \exists x_1,\dots,x_n \left(\bigwedge_i \varphi^{(r)}(x_i)\land \bigwedge_{i<j}\delta(x_i,x_j)>2r\right).
\end{align*}
In other words, there is an $r$-independent set of size $n$, with each element satisfying $\varphi^{(r)}$. 

\begin{theorem}[Gaifman]\label{thm:gaifmannormalform}
    Every first-order formula $\varphi$ is equivalent over $\mfK_\sigma$ to a Boolean combination of basic local formulae and basic local sentences. Further, the translation is effective, and if $\varphi$ is a sentence then only basic local sentences appear in the Boolean combination. 
\end{theorem}

We now proceed to the proof that order-invariant $\CFO$ contains $\FO$ on arbitrary classes of finite colored graphs. First, we show that local formulae with one free variable are expressible in $\CFO$. Second, we show that basic local sentences are equivalent to asserting the existence of independent sets in locally defined graphs, and that the existence of independent sets can be expressed in order-invariant $\CFO$.

To show that local formulae with one free variable are expressible in $\CFO$, we effectively show that we can quantify over paths of length $r$. The proof is very similar to that of {\cite[Remark 2.1]{ghasemigrange}}, so we only very briefly describe it.

\begin{lemma}\label{lem:translatinglocalformulae}
    For every $r$-local $\varphi^{(r)}(x)$ with one free variable, any valid $S$, and any $(w,i)\in S$, there exists $\varphi_{(w,i)}^{(r)}$ in $\CFO_S$ such that for any interpretation $\nu:X_S\to A$, $(A,\nu)\models \varphi_{(w,i)}^{(r)}$ if and only if $A\models \varphi^{(r)}(\nu(x_w^i))$. 
\end{lemma}

\begin{proof}
    Fix a valuation $\nu:X_S\to A$, and let $a=\nu(x_w^i)$. Each quantifier in $\varphi^{(r)}(x)$ is bounded to $B_r^A(x)$. Following {\cite[Remark 2.1]{ghasemigrange}}, each such bounded quantifier can be simulated in cluster $w$ by quantifying over paths in $A$ starting at $x_w^i$ of length $s\leq r$, with the end of the path being the quantified variable (and the set of possible endpoints is exactly $B_r^A(a)$). Because all of the variables introduced in this translation are in the cluster $w$, arbitrary edge and equality comparisons are allowed between them. Also, note that this translation is effective.

\end{proof}

The following lemma is the heart of the proof of our main containment theorem. Because we can express every $r$-local $\varphi^{(r)}(x)$, the only remaining thing to prove to show that we can express basic local formulae is to show that we can express the existence of $2r$-independent sets of size $n$ with each element satisfying $\varphi^{(r)}(x)$. The proof in general ends up being very similar to showing that there exists an independent set of size $n$, so the following lemma is the key step of the proof, and the lemma following it generalizes it suitably.

\begin{lemma}\label{lem:independensetingraph}
    For every $n$, there exists a sentence $\iota_n$ in order-invariant $\CFO$ which says that there exists an independent set of size $n$ in a finite graph $A$. 
\end{lemma}

\begin{proof}
    First, consider a finite connected graph $A$. The existence of an independent set of size $n$ in $A$ is equivalent to the existence of some $a\in A$ such that the ball $B_{2n-2}^A(a)$ contains an independent set of size $n$. Clearly the latter condition implies the former. 
    
    Now assume that $A$ contains an independent set of size $n$. If the diameter of $A$ is at most $2n-2$, then $A=B_{2n-2}^A(a)$ for every $a\in A$, so any $a$ and our assumed independent set of size $n$ witness the latter condition. If the diameter of $A$ is greater than $2n-2$, we can find a path $v_0,v_1,\dots,v_{2n-2}$ which is the shortest possible between $v_0$ and $v_{2n-2}$. Therefore, there are no edges between vertices in this path other than between successive vertices, so the set $v_0,v_2,\dots,v_{2n-2}$ is independent and contained in $B_{2n-2}^A(v_0)$. In particular, note that this second chain of reasoning does not rely on our assumed independent set of size $n$. 

    There is a $(2n-1)$-local first-order $\beta_n^{(2n-1)}(x)$ which says that $B_{2n-2}^A(x)$ contains an independent set of size $n$. Similarly, for $m<n$, we can write an $\FO$ formula $\mu_{n,m}^{(2n-1)}(x)$, which references a linear order $<$ on $A$ but is $(2n-1)$-local with respect to the edge relation on $A$, which says that
    \begin{enumerate}
        \item there is no $y$ such that $\delta^A(x,y)=2n-1$, 
        \item $x$ is the minimal element of $B_{2n-2}(x)$ with respect to $<$, and
        \item the maximum size of an independent set in $B_{2n-2}^A(x)$ is $m$.
    \end{enumerate}
    In an ordered finite graph $A$, $\mu_{n,m}^{(2n-1)}(x)$ holds exactly when $x$ is the minimal element with respect to $<$ of a connected component $C\subseteq A$ such that the maximum size of an independent set in $C$ is exactly $m<n$. The correspondence between the minimality condition and \textbf{(2)} is clear, the existence of a $y$ as described in \textbf{(1)} would imply the existence of an independent set of size $n>m$, and the non-existence of such a $y$ implies that $B_{2n-2}^A(x)$ is equal to the connected component $C$, which when combined with \textbf{(3)} completes our explanation.

    Via the same logic as \Cref{lem:translatinglocalformulae} adapted to allow for the order comparisons in \textbf{(2)}, we can write for valid $S$ and $(w,i)\in S$ a $\CFO_S$ formula $\varphi_{\mu,n,m}$ which says that $x_w^i$ satisfies $\mu_{n,m}^{(2n-1)}$. In particular, the order comparisons are possible because all variables in the local formula are in the same cluster.

    We now describe our construction of $\iota_n$. First, what we have established so far implies that $\exists x\beta_n^{(2n-1)}(x)$ is equivalent to the existence of a connected component in $A$ that contains an independent set of size $n$, and we can write a $\CFO$ sentence $\varphi_{\beta,n}$ which is equivalent to $\exists x\beta_n^{(2n-1)}(x)$. We construct $\varphi_{\beta,n}$ as
    \begin{align*}
        \varphi_{\beta_n}\equiv \exists x_\varepsilon^0\beta_{(\varepsilon,0)}^{(2n-1)}(x_{\varepsilon}^0),
    \end{align*}
    where $\beta_{(\varepsilon,0)}^{(2n-1)}(x_{\varepsilon}^0)$ is constructed as in \Cref{lem:translatinglocalformulae}. If each connected component of $A$ has no $x$ satisfying $\beta_n^{(2n-1)}(x)$, then for each connected component $C$ there is some $c\in C$ and some $m<n$ such that $\mu_{n,m}^{(2n-1)}(c)$ holds. Take $p\leq n$ and $\vec m=(m_1,\dots,m_p)$ such that $m_1+\dots+m_p\geq n$. We can write an order-invariant $\CFO$ sentence $\varphi_{\mu,\vec m,n}$ which says that there exist $x_1<\dots<x_p$ such that $\mu_{n,m_q}^{(2n-1)}(x_q)$ for every $x_q$. The construction is
    \begin{align*}
        \varphi_{\mu,\vec m,n}\equiv \exists x_\varepsilon^0 \exists x_{\alpha}^0 \dots\exists x_{\alpha^{p-1}}^0 \left(\bigwedge_{1\leq q\leq p-1}x_{\alpha^{q-1}}^0<x_{\alpha^{q}}^0\wedge \varphi_{\mu,n,m_q}(x_{\alpha^{q-1}}^0)\land \varphi_{\mu,n,m_p}(x_{\alpha^{p-1}}^0)\right).
    \end{align*}
    Finally, we have
    \begin{align*}
        \iota_n\equiv \varphi_{\beta_n}\lor \bigvee_{p,\vec m}\varphi_{\mu,\vec m,n},
    \end{align*}
    where $p$ ranges over $1\leq p\leq n$ and $\vec m$ ranges over the tuples $\vec m=(m_1,\dots,m_p)$ such that $1\leq m_i< n$ for every $i$ and $m_1+\dots+m_p\geq n$.
\end{proof}

\begin{lemma}\label{lem:localdefinableindependentsets}
    Take some $r_1$-local $\varphi^{(r_1)}(x)$ and some $r_2$-local $\psi^{(r_2)}(x,y)$. Then there exists an order-invariant $\CFO$ sentence $\iota_{\varphi,\psi,r_3,n}$ which says for every finite colored graph $A$ there is an independent set of size $n$ in the graph $G$ with domain $\{a\in A:A\models \varphi^{(r_1)}(a)\}$ and edge relation
    \begin{align*}
        \{(a,a')\in A^2:A\models (\psi^{(r_2)}(a,a')\lor \psi^{(r_2)}(a',a))\land \delta(a,a')\leq r_3\land a\neq a'\}.
    \end{align*}
\end{lemma}

\begin{proof}
    If $y$ lies in $B_{2n-1}^G(x)$, then $\delta^A(x,y)\leq r_3(2n-1)$. Also, whether $z\in A$ is a vertex of $G$ is determined by $N_{r_1}^A(z)$, so whether $(z,z')$ is an edge in $G$ is determined by $N_{\max(r_1,r_2,r_3)}^A(z,z')$. Therefore, the value of any $(2n-1)$-local formula $\phi^{(2n-1)}(x)$ in $G$ is equivalent to an $r_3(2n-1)+\max(r_1,r_2,r_3)$-local formula in $A$. Therefore the lemma follows from similar reasoning to \Cref{lem:independensetingraph}, in which the reasoning inside of each component is almost identical and the use of the ordering between the minima of components is more or less unchanged.
\end{proof}

We can now complete the proof that order-invariant $\CFO$ contains $\FO$ on colored graphs.

\begin{theorem}
    Let $\varphi$ be a first-order sentence in the vocabulary $\sigma$. Then $\varphi$ is equivalent over the colored graphs over $\sigma$ to an order-invariant $\CFO$ sentence $\psi$. In other words, on finite colored graphs, $\FO$ is contained in order-invariant $\CFO$.
\end{theorem}

\begin{proof}
    Consider a basic local sentence
    \begin{align*}
        \exists x_1,\dots,x_n \left(\bigwedge_i \alpha^{(r)}(x_i)\land \bigwedge_{i<j}\delta(x_i,x_j)>2r\right).
    \end{align*}
    Applying \Cref{lem:localdefinableindependentsets} with $\varphi^{(r_1)}=\alpha^{(r)}$, $\psi^{(r_2)}$ a tautology, and $r_3=2r$, we can write the above basic local sentence as an equivalent order-invariant $\CFO$ sentence. Thus, we can write any Boolean combination of basic local sentences as an equivalent order-invariant $\CFO$ sentence. The theorem then follows from \Cref{thm:gaifmannormalform}.
\end{proof}

The generalization to local-order-invariant $\FO$ is straightforward, which completes the proof of our main containment.

\begin{corollary}\label{cor:CFOcontainsLOinvFOcoloredgraphs}
    Let $\varphi$ be a local-order-invariant sentence over a class $\mfK\subseteq \mfF_\sigma$ of finite colored graphs in the vocabulary $\sigma\cup\{\prec\}$. Then $\varphi$ is equivalent over $\mfK$ to an order-invariant $\CFO$ sentence $\psi$. In other words, on finite colored graphs, local-order-invariant $\FO$ is contained in order-invariant $\CFO$.
\end{corollary}

\begin{proof}
    Applying \Cref{thm:gaifmannormalform} to $\varphi$, we have that $\varphi$ is a Boolean combination of basic local sentences of the form
    \begin{align*}
        \exists x_1,\dots,x_n \left(\bigwedge_i \alpha^{(r)}(x_i)\land \bigwedge_{i<j}\delta(x_i,x_j)>2r\right),
    \end{align*}
    where $\delta(x_i,x_j)>2r$ and the localization of $\alpha^{(r)}(x)$ are with respect to the Gaifman graph of $(A,\prec)$. Because $B_1^{(A,\prec)}(x)=B_2^A(x)$, then $\delta^{(A,\prec)}(x_i,x_j)>2r$ is equivalent to $\delta^A(x_i,x_j)>4r$ and $\alpha^{(r)}(x)$ is equivalent to a formula in the vocabulary $\sigma\cup \{\prec\}$ such that the domain of the quantified variables is bounded to $B_{2r}^A(x)$. Given a linear order $<$ on $A$, then for each $\alpha^{(r)}(x)$, we can replace each atom $\prec(x,y,z)$ with a formula that says that $x$ is adjacent to $y$, $x$ is adjacent to $z$, and $y<z$, via similar logic to \Cref{lem:translatinglocalformulae}. Thus, via similar logic to \Cref{lem:localdefinableindependentsets}, we have the corollary.
\end{proof}

We now conclude with several corollaries of our main containment theorem. As noted in the introduction, our first corollary generalizes \Cref{thm:aokicollapse}.

\begin{corollary}
    Let $\varphi$ be a local-order-invariant sentence over a class $\mfK\subseteq \mfF_\sigma^d$ of colored graphs in the vocabulary $\sigma\cup\{\prec\}$. Then $\varphi$ is equivalent over $\mfK$ to an $\FO$ sentence $\psi$.
\end{corollary}

\begin{definition}
    Let $\mfK\subseteq \mfF_\sigma^d$ be a class of finite structures and $\mathcal{L}$ a logic. The model checking problem for $\mathcal{L}$ on $\mfK$ is the problem of determining, with input $A\in\mfK$ and a sentence $\varphi$ in $\mathcal{L}$, whether $A\models \varphi$. We say that this model checking problem is fixed-parameter tractable if there is some computable function $f(\varphi)$ and some constant $c$ such that the model-checking problem can be solved in time $f(\varphi)\cdot O(|A|^c)$.  
\end{definition}

\begin{theorem}[Ghasemi--Grange]\label{thm:CFOFPT}
    The model checking problem for order-invariant $\CFO$ on a class $\mfK\subseteq \mfF_\sigma^d$ is fixed-parameter tractable.
\end{theorem}

\begin{corollary}\label{cor:localorderFPTcoloredgraphs}
    The model checking problem for local-order-invariant $\FO$ on a class $\mfK\subseteq \mfF_\sigma^d$ of colored graphs is fixed-parameter tractable.
\end{corollary}

\begin{proof}
    Our construction of a $\CFO$ sentence expressing a basic local sentence in the vocabulary $\sigma\cup\{\prec\}$ is effective. Also, as is noted in \Cref{thm:gaifmannormalform}, the translation from an arbitrary local-order-invariant $\FO$ sentence to Gaifman normal form is effective. Combining these, we have an effective translation from local-order-invariant $\FO$ sentences to order-invariant $\CFO$ sentences, so \Cref{thm:CFOFPT} extends to local-order-invariant $\FO$. 
\end{proof}

\section{Generalizations via Interpretations}\label{sec:interpretations}

Our main containment result was proved for colored graphs, because this is the setting in which $\CFO$ was originally defined and in which \Cref{thm:CFOcollapseoncoloredgraphs} applies. We now generalize \Cref{thm:CFOcollapseoncoloredgraphs}, \Cref{cor:CFOcontainsLOinvFOcoloredgraphs}, and \Cref{cor:localorderFPTcoloredgraphs} to general vocabularies, affirming \Cref{conj:CFOcontainsFO}.

We begin by defining interpretations themselves.

\begin{definition}
    A $d$-dimensional interpretation $I$ between $\mfK\subseteq \mfK_\sigma$ and $\mfK_\tau$ is a domain formula $\delta(\vec x)$ with $\vec x$ of length $d$, an equivalence formula $\varepsilon(\vec x,\vec y)$ with $\vec x,\vec y$ of length $d$, and a relation formula $\rho_R(\vec x_1,\dots,\vec x_{a_R})$ for each $R\in\tau$ with each $\vec x_i$ of length $d$. Each of these formulae is in the vocabulary $\sigma$. For a given $A\in \mfK$, we write $D_I(A)$ for the set $\{\vec a\in A^d:\delta(\vec a)\}$, and we require $\varepsilon$ to define an equivalence relation on $D_I(A)$ for every $A\in \mfK$, and further we require that if we have $\vec a_1,\dots,\vec a_{a_R}\in D_I(A)$ and $\vec a_1',\dots,\vec a_{a_R}'\in D_I(A)$ such that $A\models \varepsilon(\vec a_i,\vec a_i')$ for every $i$, then $A\models \rho_R(\vec a_1,\dots,\vec a_{a_R})$ if and only if $A\models \rho_R(\vec a_1',\dots,\vec a_{a_R}')$. We define $I(A)$ as the structure over $\tau$ with underlying set $D_I(A)/\varepsilon$, with relations $R^{I(A)}([\vec a_1],\dots,[\vec a_{a_R}])$ such that $R^{I(A)}([\vec a_1],\dots,[\vec a_{a_R}])$ in $I(A)$ if and only if $A\models \rho_R(\vec a_1,\dots,\vec a_{a_R})$. 
\end{definition}

\begin{theorem}\label{thm:allcollapsesonallboundeddegree}
    \begin{enumerate}
        \item Let $\sigma$ be an arbitrary vocabulary, and let $\mfK\subseteq \mfF_\sigma$. Every sentence $\varphi$ which is local-order-invariant over $\mfK$ is equivalent on $\mfK$ to an order-invariant $\CFO$ sentence $\psi$.
        \item Conversely, let $\mfK\subseteq \mfF_\sigma^d$. Every $\CFO$ sentence $\varphi$ which is order-invariant over $\mfK$ is equivalent over $\mfK$ to a first-order sentence $\psi$.
    \end{enumerate}
     In particular, order-invariant $\CFO$, local-order-invariant $\FO$, and plain $\FO$ are equivalent over any class of finite structures of bounded degree.
\end{theorem}

\begin{proof}
    The containment of local-order-invariant $\FO$ in order-invariant $\CFO$ on classes of finite structures over arbitrary vocabularies is an almost verbatim generalization of the proof of \Cref{cor:CFOcontainsLOinvFOcoloredgraphs}. 
    
    We now describe the generalization of the collapse result on classes of bounded degree. For structures with one element the translation is trivial, so we assume from here that our structures have at least two elements. Take a structure $A$ over $\sigma$. Define the colored simple graph $\Inc(A)$ as follows.
    \begin{enumerate}
        \item There is a color $U$ for the domain of $A$, a relation color $T_R$ for each $R\in \sigma$, and a position color $P_{R,i}$ for each $1\leq i\leq a_R$ where $a_R$ is the arity of $R$.
        \item For each element $a$ of $A$, there is a vertex $u_a$ of $\Inc(A)$ colored $U$.
        \item For each relation $R$ and tuple $\vec a$ of elements of $A$ such that $A\models R(\vec a)$, there is a vertex $t_{R,\vec a}$ of $\Inc(A)$ colored $T_R$.
        \item For each vertex $t_{R,\vec a}$, and each $1\leq i\leq a_R$, there is a vertex $p_{R,\vec a,i}$ of $\Inc(A)$ colored $P_{R,i}$. 
        \item For every vertex $t_{R,\vec a}$ and $1\leq i\leq a_R$, we have $E(t_{R,\vec a},p_{R,\vec a,i})$ and $E(p_{R,\vec a,i},u_{a_i})$ in $\Inc(A)$.
    \end{enumerate}
    There is a uniform first-order interpretation $I$ such that $I(\Inc(A))\cong A$. This construction $\Inc$ preserves boundedness of degree, in the sense that for all $\mfK\subseteq \mfF_\sigma^d$, the class $\Inc(\mfK)=\{\Inc(A):A\in\mfK\}$ is of bounded degree. To see why, take $a_\sigma=\max\{a_R:R\in \sigma\}$. Then for all $a\in A$ and any $R\in \sigma$, there are at most $(d+1)^{a_\sigma}$ tuples $\vec a$ of elements of $A$ such that $R(\vec a)$ for some $R\in\sigma$. Therefore each $u_a$ has degree at most
    \begin{align*}
        \sum_{R\in\sigma}a_R(d+1)^{a_R}\leq |\sigma|a_\sigma(d+1)^{a_\sigma},
    \end{align*}
    each $t_{R,\vec a}$ has degree exactly $a_R$, and each $p_{R,\vec a,i}$ has degree exactly $2$, so $\Inc(A)$ has degree at most
    \begin{align*}
        \max(|\sigma|a_\sigma(d+1)^{a_\sigma},2).
    \end{align*}
    Every $\varphi_1$ which is an order-invariant $\CFO$ sentence on $\mfK\subseteq \mfF_\sigma^d$ can be translated to an order-invariant $\CFO$ sentence $\varphi_2$ on $\Inc(\mfK)$. Such a translation can proceed inductively, by encoding each atom $R(\vec x)$ via similar logic to \Cref{lem:translatinglocalformulae}, translating guarded quantification to quantification over element vertices of distance exactly $4$ again by similar logic to \Cref{lem:translatinglocalformulae}, and preserving order comparisons and unguarded quantification. By \Cref{thm:CFOcollapseoncoloredgraphs}, this translated sentence is equivalent on $\Inc(\mfK)$ to a first-order sentence $\varphi_3$. Finally, we claim that $\Inc$ is a uniform first-order interpretation, so $\varphi_3$ can be translated to a first-order sentence $\varphi_4$ which is equivalent on $\mfK$ to $\varphi_1$. We now describe this interpretation explicitly.

    Let
    \begin{align*}
        \mathcal S=\{U\}\cup \{T_R:R\in \sigma\}\cup\{P_{R,i}:R\in \sigma,i\leq a_R\}.
    \end{align*}
    The equality type of a tuple $\vec a$ of elements of $A$ is the set $\{(i,j):a_i=a_j\}$. Choose $h$ to be large enough so that there are at least $|\mathcal S|$ equality types for elements of $\{0,1\}^h$. For each $s\in \mathcal S$, choose the equality type $\tau_s$ of some element of $\{0,1\}^h$, and construct a formula $\eta_{\tau_s}(\vec y)$ which says that $\vec y$ has equality type $\tau_s$. We define an $(a_\sigma+h)$-dimensional interpretation $J$ between $\mfK$ and $\mfK_{\{E\}\cup \mathcal S}$ such that $J(A)\cong \Inc(A)$ for every $A\in \mfK$ of size at least two. Informally, the first $a_\sigma$ coordinates will represent the elements themselves, and the latter $h$ coordinates will represent the sort from $\mathcal S$. More precisely, we have the following.
    \begin{enumerate}
        \item If $\eta_{\tau_s}(\vec y)$ holds with $s=U$, then $\delta(\vec x,\vec y)$ is true, and $\varepsilon(\vec x,\vec y,\vec x',\vec y')$ holds if and only if $x_1=x_1'$ and $\eta_{\tau_s}(\vec y')$ holds.
        \item If $\eta_{\tau_s}(\vec y)$ holds with $s=T_R$, then $\delta(\vec x,\vec y)$ holds if and only if $R(x_1,\dots,x_{a_R})$ holds, and $\varepsilon(\vec x,\vec y,\vec x',\vec y')$ holds if and only if $x_i=x_i'$ for $i\leq a_R$ and $\eta_{\tau_s}(\vec y')$ holds.
        \item If $\eta_{\tau_s}(\vec y)$ holds with $s=P_{R,i}$, then $\delta(\vec x,\vec y)$ holds if and only if $R(x_1,\dots,x_{a_R})$ holds, and $\varepsilon(\vec x,\vec y,\vec x',\vec y')$ holds if and only if $x_i=x_i'$ for $i\leq a_R$ and $\eta_{\tau_s}(\vec y')$ holds.
        \item If $\eta_{\tau_s}(\vec y)$ does not hold for any $s\in \mathcal S$, then $\delta(\vec x,\vec y)$ does not hold.
        \item For $s\in \mathcal S$, the formula $\rho_s(\vec x,\vec y)$ holds if and only if $\eta_{\tau_s}(\vec y)$ holds.
        \item The formula $\rho_E(\vec x,\vec y,\vec x',\vec y')$ holds if and only if either
        \begin{enumerate}
            \item $x_1=x_i'$ (or $x_i=x_1'$), $\eta_{\tau_s}(\vec y)$ and $\eta_{\tau_{s'}}(\vec y')$ hold with $s=U$ and $s'=P_{R,i}$ (or vice versa); or
            \item $x_i=x_i'$ for $i\leq a_R$, and $\eta_{\tau_s}(\vec y)$ and $\eta_{\tau_{s'}}(\vec y')$ hold with $s=T_R$ and $s'=P_{R,i}$ (or vice versa).
        \end{enumerate}
    \end{enumerate}
    That $J(A)\cong \Inc(A)$ is clear, so we are done.
\end{proof}

\begin{corollary}\label{cor:localorderFPT}
    For any vocabulary $\sigma$, the model checking problem for local-order-invariant $\FO$ on a class $\mfK\subseteq \mfF_\sigma^d$ is fixed-parameter tractable.
\end{corollary}

\begin{proof}
    The translation from a local-order-invariant $\FO$ sentence to an order-invariant $\CFO$ sentence is effective, and the translation from $\varphi_1$ to $\varphi_2$ in \Cref{thm:allcollapsesonallboundeddegree} is effective, so the corollary follows from \Cref{thm:CFOFPT}.
\end{proof}

A question of interest is whether local-order-invariant $\FO$ is as expressive as order-invariant $\FO$, which is nontrivial because they have the same expressiveness on classes with Gaifman graphs of bounded diameter by a result of \cite{aoki2025}. In particular, an answer to this question would answer the important question of whether order-invariant $\FO$ collapses to $\FO$ on classes of bounded degree. We give an equivalent formulation for this question in terms of interpretations.

\begin{definition}
    We say that a logic $\mathcal{L}$ is closed under first-order interpretations if, for every interpretation $I$ from $\mfK\subseteq \mfK_\sigma$ to $\mfK_\tau$ and every $\mathcal{L}$-definable $\mfJ\subseteq I(\mfK)$, the preimage $I^{-1}(\mfJ)$ is $\mathcal{L}$-definable in $\mfK$.
\end{definition}

We defer the proof of the following proposition to \Cref{sec:proofsforsec4}.

\begin{proposition}\label{prop:interpretationformulationofequality}
    Local-order-invariant $\FO$ is closed under first-order interpretations on all classes of finite structures if and only if it has the same expressive power as order-invariant $\FO$ on all classes of finite structures
\end{proposition}

The following corollary follows from similar techniques. We defer the definition of tree-width and path-width to \Cref{sec:proofsforsec4}, as well as the proof of the corollary.

\begin{corollary}\label{cor:collapseonboundedtreewidth}
    If local-order-invariant $\FO$ has the same expressiveness as $\FO$ on every class of bounded tree-width (resp. bounded path-width), then order-invariant $\FO$ has the same expressiveness as $\FO$ on every class of bounded tree-width (resp. path-width).
\end{corollary}

One of the primary reasons that we believe this analysis of invariant logics via interpretations has been under-developed is that it yields straightforwardly a collapse of order-invariant $\FO$ to $\FO$ when restricted to classes of bounded shrub-depth. The proof is simply a combination of results of \cite{eickmeyerElberfeldHarwath2017}, \cite{chenFlum2020}, and \cite{ohlmannPilipczukPrzybyszewskiTorunczyk2023}.

\begin{definition}
    A graph $G$ has tree-depth at most $d$ if there is a rooted forest $F$ of height $d$ such that every edge of $G$ joins an ancestor and a descendant in $F$. A class of graphs has bounded tree-depth if there is some $d$ such that each element has tree-depth at most $d$. A class $\mfK$ of structures over an arbitrary vocabulary has bounded tree-depth if the class of graphs $\{\mcG(A):A\in \mfK\}$ is of bounded tree-depth. 
\end{definition}

\begin{definition}
    A tree-model of $m$ colors and depth $d$ for a graph $G$ is a pair $(T,S)$ of a rooted tree $T$ of height $d$ and a set $S\subseteq \{1,\dots,m\}^2\times \{1,\dots,d\}$ such that
    \begin{enumerate}
        \item the length of each root to leaf path in $T$ is exactly $d$,
        \item the set of leaves of $T$ is exactly the set of vertices of $G$,
        \item each leaf of $T$ is assigned one of the $m$ colors, and
        \item for all $i,j,l$, we have $(i,j,l)\in S$ if and only if $(j,i,l)\in S$, and 
        \item for any two vertices $u,v$ of $G$ with colors $i,j$ respectively such that the distance between $u$ and $v$ in $T$ is $2l$, there is an edge between $u$ and $v$ in $G$ if and only if $(i,j,l)\in S$.
    \end{enumerate}
    A class of graphs is of bounded shrub-depth if and only if there exist $m,d$ such that every graph in the class has a tree-model of $m$ colors and depth $d$.
\end{definition}

The following results are from \cite{ohlmannPilipczukPrzybyszewskiTorunczyk2023,eickmeyerElberfeldHarwath2017,chenFlum2020}, though it is worth noting that \Cref{thm:ChenFlum} was originally claimed in \cite{gajarskyHlineny2015}.

\begin{definition}
    A simple interpretation is a $1$-dimensional interpretation with $\varepsilon(x,y)\equiv(x=y)$. A transduction $T$ consists of an integer $k$, a binary relation $S$, unary relations $U_1,\dots,U_n$, and a simple interpretation $I_T$ from $\mfK\subseteq \mfK_{\sigma\cup\{S,U_1,\dots,U_n\}}$ to $\mfK_\tau$. Given $A\in \mfK$, write $A_k$ for the set $\{1,\dots,k\}\times A$, and write $(A_k,S)$ for the structure over $\sigma\cup\{S\}$ such that as structures over $\sigma$ each $\{i\}\times A$ is isomorphic to $A$, and $S$ is the equivalence relation on $A_k$ such that $S((i,a),(j,a'))$ if and only if $a=a'$. We define the output $T(A)$ to be the set of structures of the form $I_T(A_k,S,U_1,\dots,U_n)$ where each $U_i$ is some unary relation on $(A_k,S)$.
\end{definition}

\begin{theorem}[Ohlmann--Pilipczuk--Przybyszewski--Toru\'nczyk]\label{thm:OPPT}
    Let $\mfK$ be a class of finite graphs of bounded shrub-depth. There exist
    \begin{enumerate}
        \item a finite relational vocabulary $\sigma$ of arity at most 2,
        \item an integer $d$,
        \item a transduction $T$ from $\mfK$ to $\mfF_\sigma$, 
        \item a unary $P$ in $\sigma$, and
        \item a first-order interpretation $I$ from $\mfF_\sigma$ to the class of finite graphs,
    \end{enumerate}
    such that, for every $G\in \mfK$, some output $B\in T(G)$ has tree-depth at most $d$, the elements in $P$ are the original vertices of $G$, and $D_I(B)$ is $\{b\in B:B\models P(b)\}$.
\end{theorem}

\begin{theorem}[Eickmeyer--Elberfeld--Harwath]\label{thm:EEH}
    Let $\varphi$ be an order-invariant sentence on $\mfF_\sigma$, and let $\mfF_\sigma^{\mathrm{td},d}$ be the class of finite $\sigma$ structures of tree-depth at most $d$. Then there exists a first-order sentence $\psi$ which is equivalent to $\varphi$ on $\mfF_\sigma^{\mathrm{td},d}$.
\end{theorem}

\begin{theorem}[Chen--Flum]\label{thm:ChenFlum}
    Let $\mfK$ be a class of finite graphs of bounded shrub-depth, and let $\varphi$ be an $\MSO$ sentence. Then $\varphi$ is equivalent on $\mfK$ to some first-order sentence $\psi$.
\end{theorem}

\begin{theorem}\label{thm:boundedshrubdepthcollapse}
    Let $\mfK$ be a class of finite graphs of bounded shrub-depth, and let $\varphi$ be order-invariant over the class of finite graphs. Then $\varphi$ is equivalent on $\mfK$ to a first-order sentence $\psi$.
\end{theorem}

\begin{proof}
    The proof is a relatively straightforward combination of \Cref{thm:EEH}, \Cref{thm:OPPT}, \Cref{thm:ChenFlum}, the closure of $\MSO$ under interpretations, and the definability of tree-depth in $\MSO$. The full proof can be found in \Cref{sec:proofsforsec4}.
\end{proof}

To attempt to generalize \Cref{thm:boundedshrubdepthcollapse} past graphs, one might attempt to apply similar methods to \Cref{thm:allcollapsesonallboundeddegree}. A first step might be to generalize \Cref{thm:boundedshrubdepthcollapse} to colored graphs. 
    
Indeed, let $\mfK$ be a class of finite colored graphs of bounded shrub-depth, and let $\varphi$ be order-invariant over the class of finite colored graphs. Then $\varphi$ is equivalent on $\mfK$ to a first-order sentence $\psi$. The proof of the above result generalizes to colored graphs straightforwardly. In slightly more detail, we can add finitely many unary color relations as unary predicates, track them through the transduction and interpretation from \Cref{thm:OPPT}, apply \Cref{thm:EEH} identically while taking these unary relations into account, and apply \Cref{thm:ChenFlum} identically while taking these unary relations into account.

\begin{remark}\label{rem:nouniformencoding}
    However, the proof of \Cref{thm:allcollapsesonallboundeddegree} cannot generalize to classes of bounded shrub-depth. That is, there is no uniform encoding of finite relational structures into finitely colored graphs such that 
    \begin{enumerate}
        \item $A$ is first-order interpretable in its encoding, and
        \item if $\mfK$ is a class of bounded shrub-depth then the encoding takes $\mfK$ to a class of colored graphs of bounded shrub-depth.
    \end{enumerate}
    In particular, there is an obstruction even for finite linear orders.
\end{remark}

\section{Deterministic Frozen-Witness Logic}\label{sec:frozenwitness}

We conclude with an application of our collapse of local-order-invariant $\FO$ to $\FO$ on classes of bounded degree. This recovers, on classes of bounded degree, the collapse of deterministic frozen-witness $\FO$ to $\FO$ conjectured by Abiteboul, Simon, and Vianu \cite{abiteboulSimonVianu1990} (which was disproved on arbitrary finite structures by Otto \cite{otto2000}). The proof takes inspiration from the erroneous proof of \cite{aoki2025} that epsilon-invariant $\FO$ is contained in local-order-invariant $\FO$ on classes of bounded degree. However, we apply the theorem to a weaker logic (deterministic frozen-witness logic), make substantial repairs to the proof, and remove the bounded degree hypothesis. Further, the erroneously claimed collapse of epsilon-invariant $\FO$ to $\FO$ on classes of bounded degree in \cite{aoki2025} inherits the necessity of invariance on an ambient substructure-closed class from \Cref{thm:aokicollapse}, whereas because we have removed that hypothesis in \Cref{thm:allcollapsesonallboundeddegree}, our collapse inherits the lack of the hypothesis. 

We begin by introducing deterministic frozen-witness logic. Our presentation is slightly non-standard, but it is clearly equivalent to the standard definition, is close to the presentation of \cite{otto2000}, and is much more convenient for the purposes of the following proofs.

\begin{definition}
    Let $\varphi(\vec x,y)$ be a first-order formula in the vocabulary $\sigma$, and let $A\in \mfK_\sigma$. We say that a relation $R(\vec x,y)$ on $A$ is a witness predicate for $\varphi$ if it satisfies the following conditions.
    \begin{enumerate}
        \item For every tuple $\vec a$ of elements of $A$ and $a'\in A$, then $R(\vec a,a')$ implies $A\models \varphi(\vec a,a')$.
        \item If there exists $a'\in A$ such that $A\models \varphi(\vec a,a')$, then there is a unique $a''\in A$ such that $R(\vec a,a'')$.
        \item If there does not exist $a'\in A$ such that $A\models \varphi(\vec a,a')$, then there is no $a''\in A$ such that $R(\vec a,a'')$.
    \end{enumerate}
\end{definition}

\begin{definition}\label{def:frozenwitnessprogram}
    A frozen witness program over a vocabulary $\sigma$ is a finite sequence of relations $R_1,\dots,R_n$ with each $R_i$ not in $\sigma$, a sequence of first-order formulae $\varphi_1,\dots,\varphi_n$ with $\varphi_i$ a formula in the vocabulary $\sigma\cup\{R_1,\dots,R_{i-1}\}$, and a first-order sentence $\psi$ in the vocabulary $\sigma\cup\{R_1,\dots,R_n\}$. A run of the program on a structure $A\in \mfK_\sigma$ consists of interpreting $R_1,\dots,R_n$ in order, such that the interpretation of $R_i$ is a witness predicate for $\varphi_i$ on $(A,R_1,\dots,R_{i-1})$. Finally, the program evaluates $\psi$ once it has fixed all of these interpretations.

    We say that a frozen witness program is deterministic on $\mfK\subseteq \mfK_\sigma$ if, for every $A\in \mfK$, all runs of the program evaluate $\psi$ with the same truth value. If every run evaluates $\psi$ to be true, then we write $A\models \psi$.
\end{definition}

The following lemma is from \cite{aoki2025}.

\begin{lemma}[Aoki]\label{lem:neighborhoodorderdefinition}
    For every $r$, there exists a formula $\lambda_r(x,y,z)$ in the vocabulary $\sigma\cup\{\prec\}$ which linearly orders $B_r^A(x)$ when $x$ is fixed and when $\prec$ is a local order.
\end{lemma}


\begin{definition}
    We write $\tp_k(A)$ for the set of first-order sentences $\varphi$ of quantifier rank at most $k$ such that $A\models\varphi$. Given a tuple $\vec a$ of elements of $A$ and $r\in\omega$, we write $\tp_k^r(A,\vec a)$ for the set of first-order $r$-local formulae $\varphi^{(r)}(\vec x)$ of quantifier rank at most $k$ such that $A\models \varphi^{(r)}(\vec a)$.
\end{definition}

The following is a well-known folklore result {\cite[Theorem 3.15]{libkin2004}}.

\begin{theorem}\label{thm:finitelymanytypes}
    For every $r,k,n\in\omega$, the sets $\{\tp_k(A):A\in\mfK_\sigma\}$ and
    \begin{align*}
        \{\tp_k^r(A,\vec a):A\in\mfK_\sigma\text{ and }\vec a\in A^n\}
    \end{align*}
    are both finite.
\end{theorem}

The following two lemmas form the core of the proof of our containment of deterministic frozen-witness $\FO$ in local-order-invariant $\FO$. The basic idea is that, given a formula $\varphi(\vec x,y)$ for which we wish to define a witness predicate, we can choose a collection $D$ of witnesses of bounded size, where the size depends on only the formula, such that for any $\vec x$ there is a witness $y$ which is distance at most $r$ from some point in either $D$ or $\vec x$. Then, we use a local order to linearly order the elements distance at most $r$ from some point in $D$ or $\vec x$, and choose the minimum valid witness using this linear order. We establish the existence of $D$ in the following lemma, and we establish this choice via local order in the lemma after.

\begin{lemma}\label{lem:goodtuples}
    Let $\theta(\vec x,y,\vec c)$ be a first-order formula in some vocabulary $\tau$. Then there are integers $r,n$, and a first-order formula $\gamma_\theta(\vec c,\vec d)$ with $\vec d$ of length $n$ such that the following hold for every $A\in \mfK_\tau$.
    \begin{enumerate}
        \item For every tuple $\vec a$ of elements of $A$ of the same length as $\vec c$, there is some tuple $\vec a'$ such that $A\models \gamma_\theta(\vec a,\vec a')$.
        \item If $A\models \gamma_\theta(\vec a,\vec a')$, $\vec b$ is a tuple of elements of $A$, and there exists $b'\in A$ such that $\theta(\vec b,b',\vec a)$, then $A\models \theta(\vec b,b'',\vec a)$ for some $b''\in B_r^A(\vec a,\vec b)$ or for some $b''\in B_{2r}^A(\vec a')$.
    \end{enumerate}
\end{lemma}

\begin{proof}
    Put $\theta(\vec x,y,\vec c)$ in Gaifman normal form, let the radius be $s$, let $r=2s$, and let the maximal quantifier rank of the $s$-local formulae in the Gaifman normal form be $k$. By \Cref{thm:finitelymanytypes}, we can enumerate the quantifier rank $k$ $s$-local types $\tau_1,\dots,\tau_m$. Let $\ell_{\vec x}$ be the length of $\vec x$, let $\ell_{\vec c}$ be the length of ${\vec c}$. Let $n=m(\ell_{\vec x}+\ell_{\vec c}+1)$. 
    
    Denote by $D$ the set $\{d_1,\dots,d_n\}$. Then $\gamma_\theta(\vec c, \vec d)$ says that for each $\tau_i$, either
    \begin{enumerate}
        \item $D$ contains a $2r$-independent subset of size $\ell_{\vec x}+\ell_{\vec c}+1$ such that $\tp_k^s(A,d)=\tau_i$ for every $d$ in this subset, or
        \item $D$ contains a $2r$-independent set of size at most $\ell_{\vec x}+\ell_{\vec c}$, and this subset is maximal among the $2r$-independent subsets of $\{a\in A:\tp_k^s(A,a)=\tau_i\}$.
    \end{enumerate}
    Note that this is first-order expressible because $D$ is of bounded size. 

    Now suppose that $A\models \theta(\vec b,b',\vec a)$ for tuples $\vec b,\vec a$ of elements of $A$ and $b'\in A$, and suppose that we have constructed $\vec a'$ such that $A\models\gamma_\theta(\vec a,\vec a')$. If $b'$ is in $B_r^A(\vec a,\vec a',\vec b)$, then we are done. If $b'$ is not in $B_r^A(\vec a,\vec a',\vec b)$, then if \textbf{(1)} from our earlier dichotomy for $D$ (where $D$ is $\{a'_1,\dots,a'_n\}$) holds for $\tau_i=\tp_k^s(A,b')$, there is some $d\in D$ which is distance greater than $r$ from every element of $\vec b$ and $\vec a$. Therefore, $\tp_k^s(A,\vec b,\vec a,b')=\tp_k^s(A,\vec b,\vec a,d)$, so $A\models \theta(\vec b,b',\vec a)$ if and only if $A\models \theta(\vec b,d,\vec a)$, and thus $A\models \theta(\vec b,d,\vec a)$ with $d\in B_{2r}(\vec a,\vec a',\vec b)$ because it is an element of $\vec a'$. If instead \textbf{(2)} from the dichotomy holds, then every $b''\in A$ such that $\tp_k^s(A,b')=\tp_k^s(A,b'')$ is in $B_{2r}^A(\vec a')$ by maximality of the relevant $2r$-independent subset of $D$. 
\end{proof}

\begin{lemma}\label{lem:localorderselectordefinition}
    Let $\theta(\vec x,y,\vec c)$ be a first-order formula in the vocabulary $\sigma\cup\{\prec\}$, and construct $r$, $n$, and $\gamma_\theta$ as in \Cref{lem:goodtuples}. Then there is a formula $\chi_\theta(\vec x,y,\vec c,\vec d)$ such that for any $A\in\mfF_\sigma$, any local order $\prec$ on $A$, and any $\vec a,\vec a'$ such that $(A,\prec)\models \gamma_\theta(\vec a,\vec a')$, the relation $R_{\vec a,\vec a'}(\vec x,y)\equiv \chi_\theta(\vec x,y,\vec a,\vec a')$ is a witness predicate for $\theta(\vec x,y,\vec a)$.
\end{lemma}

\begin{proof}
    Construct $\chi_\theta$ as follows. 
    \begin{enumerate}
        \item If there are no $y$ such that $\theta(\vec x,y,\vec c)$, then set $\chi_\theta(\vec x,y\vec c,\vec d)$ to be false for every $y$.
        \item If some $d_i$ in $\vec d$ satisfies $\theta(\vec x,d_i,\vec c)$, then let $j$ be the least $i$ such that $\theta(\vec x,d_i,\vec c)$, and set $\chi(\vec x,y,\vec c,\vec d)$ to be true if and only if $y=d_j$
        \item If there is some $y$ such that $\theta(\vec x,y,\vec c)$, but no $d_i$ satisfies $\theta(\vec x,d_i,\vec c)$, then there exists some $y$ in $B_r^{(A,\prec)}(\vec x,\vec c)\cup B_{2r}^{(A,\prec)}(\vec d)$ such that $\theta(\vec x,y,\vec c)$. 
        
        We construct a linear order on $B_{2r}^{(A,\prec)}(\vec x,\vec c,\vec d)$ as follows. Define the tuple $\vec z$ as the concatenation 
        \begin{align*}
            z_1,\dots,z_k=x_1,\dots,x_m,c_1,\dots,c_l,d_1,\dots,d_n.
        \end{align*}
        For each $p\in B_{2r}^{(A,\prec)}(\vec x,\vec c,\vec d)$, let $\kappa(p)$ be the least $i$ such that $p\in B_{2r}^{(A,\prec)}(z_i)$. For $p,q\in B_{2r}^{(A,\prec)}(\vec x,\vec c,\vec d)$, let $p$ precede $q$ if and only if $\kappa(p)<\kappa(q)$ or $\kappa(p)=\kappa(q)$ and $\lambda_{4r}(z_i,p,q)$ holds, with $\lambda_{4r}$ as in \Cref{lem:neighborhoodorderdefinition}. This works because $B_{2r}^{(A,\prec)}(\vec a)$, for any $A\in \mfK_\sigma$, any tuple $\vec a$ of elements of $A$, any local order $\prec$ on $A$, and any $r\in \omega$, is equal to $B_{4r}^A(\vec a)$.
        
        Now, using this linear order, which is definable in the vocabulary $\sigma\cup\{\prec\}$, simply take the least element $y$ of $B_r^{(A,\prec)}(\vec x,\vec c)\cup B_{2r}^{(A,\prec)}(\vec d)$ such that $\theta(\vec x,y\vec c)$.
    \end{enumerate}
    These cases are exhaustive, and so we are done.
\end{proof}

We now complete the proof of our main containment theorem for this section and give several corollaries.

\begin{theorem}\label{thm:LOinvFOcontainsfrozenwitness}
    Take some deterministic frozen witness program on $\mfK\subseteq \mfF_\sigma$ consisting of $R_1,\dots,R_n$, $\varphi_1,\dots,\varphi_n$, and $\psi$. Then there exists some local-order-invariant $\FO$ sentence $\phi$ on $\mfK$ such that for every $A\in\mfK$, $A\models \psi$ if and only if $A\models \phi$. 
\end{theorem}

\begin{proof}
    We construct, by induction on $i$, formulae $\rho_i(\vec x_i,y_i,\vec d_{\leq i})$ in the vocabulary $\sigma\cup\{\prec\}$, where $\vec d_{\leq i}=(\vec d_1,\dots,\vec d_i)$, such that for any suitable $\vec d_{\leq i}$ (in a sense which we will define), the relations $R_{\rho_i}(\vec x,y)\equiv \rho_i(\vec x,y,\vec d_{\leq i})$ are suitable interpretations for the $R_i$ in the program.

    For $i=1$, we simply apply \Cref{lem:localorderselectordefinition} to $\varphi_1(\vec x_1,y_1)$. This gives us a formula $\gamma_1(\vec d_1)$ and a formula $\rho_1(\vec x_1,y_1,\vec d_1)$ defining a witness predicate for $\varphi_1$. 
    
    Now assume that $\rho_1,\dots,\rho_{i-1}$ have been constructed. In the formula $\varphi_i$, replace each $R_j(\vec x,y)$ with $\rho_j(\vec x,y,\vec d_{\leq j})$. Call this new formula $\varphi'_i(\vec x_i,y_i,\vec d_{\leq i-1})$. Then using \Cref{lem:localorderselectordefinition} again, we can construct $\gamma_i(\vec d_{\leq i})$ and $\rho_i(\vec x_i,y,\vec d_{\leq i})$ which define a witness predicate for $\varphi_i$.

    Now construct a formula $\psi'(\vec d_{\leq n})$ by replacing each $R_i(\vec x,y)$ in $\psi$ with $\rho_i(\vec x,y,\vec d_{\leq i})$. Finally, construct the sentence
    \begin{align*}
        \psi''\equiv \exists\vec d_1\dots\exists \vec d_n\left(\psi'(\vec d_{\leq n})\land \bigwedge_{i\leq n}\gamma_i(\vec d_{\leq i})\right).
    \end{align*}
    For every $A\in\mfK$ and every local order $\prec$ on $A$, there exists a tuple $\vec d_{\leq n}$ satisfying the conditions $\gamma_i(\vec d_{\leq i})$ for each $i$, and for any such tuple, the $\rho_i(\vec x_i,y_i,\vec d_{\leq i})$ define witness predicates, so for any local order $\prec$ on $A$, $\psi''$ says that $A\models \psi$ on a certain run of the program. Therefore, because $\psi$ evaluates on $\mfK$ in the same way regardless of the chosen witness predicates, $A\models\psi$ if and only if $A\models \psi''$ for any local order $\prec$, so $\psi''$ is equivalent to $\psi$ and is local-order-invariant on $\mfK$.
\end{proof}

\begin{corollary}
    Take some deterministic frozen witness program on $\mfK\subseteq \mfF_\sigma^d$ consisting of $R_1,\dots,R_n$, $\varphi_1,\dots,\varphi_n$, and $\psi$. Then there exists some first-order sentence $\phi$ on $\mfK$ such that for every $A\in\mfK$, $A\models \psi$ if and only if $A\models \phi$.
\end{corollary}

\begin{proof}
    This is simply a combination of \Cref{thm:allcollapsesonallboundeddegree} and \Cref{thm:LOinvFOcontainsfrozenwitness}.
\end{proof}

\begin{corollary}
    For any $\mfK\subseteq \mfF_\sigma^d$, the model-checking problem for deterministic frozen-witness programs is fixed-parameter tractable on $\mfK$.
\end{corollary}

\begin{proof}
    Our translation in \Cref{thm:LOinvFOcontainsfrozenwitness} is effective, in particular because the translation to Gaifman normal form is effective by \Cref{thm:gaifmannormalform}, so this is simply an application of \Cref{cor:localorderFPT}.
\end{proof}

\section{Conclusion}

We have shown that order-invariant $\CFO$ contains local-order-invariant $\FO$, affirming \Cref{conj:CFOcontainsFO}. We have also shown that local-order-invariant $\FO$ contains deterministic frozen-witness $\FO$, salvaging the proof strategy of a claimed result of \cite{aoki2025}. Together, these prove the equality of $\FO$, deterministic frozen-witness $\FO$, local-order-invariant $\FO$, and order-invariant $\CFO$ on classes of bounded degree. We also used our method of interpretations to extend the results beyond colored graphs and prove the collapse of order-invariant $\FO$ when restricted to classes of bounded shrub-depth. There are several natural remaining questions, such as whether local-order-invariant $\FO$ and order-invariant $\FO$ have the same expressiveness, or whether interpretations can be used in other ways to study invariant logics.

\bibliography{bibliography}

\appendix

\section{Proofs for \Cref{sec:interpretations}}\label{sec:proofsforsec4}

\begin{proof}[Proof of \Cref{prop:interpretationformulationofequality}]
    First, we note that order-invariant $\FO$ is closed under interpretations. Take a $d$-dimensional interpretation $I$ between $\mfK\subseteq \mfF_\sigma$ and $\mfF_\tau$ defined by $\delta$, $\varepsilon$, and $\rho_R$ for each $R\in \tau$. Given $A\in \mfK$ and a linear order $<$ on $A$, we can define a linear order on $I(A)$ by first taking the lexicographic order on $A^d$ induced by $<$ and then ordering equivalence classes of $D_I(A)/\varepsilon$ by their minimal elements. That is, we will have $[\vec a_1]<_{I(A)}[\vec a_2]$ if $\vec a_1'$ precedes $\vec a_2'$ in the lexicographic order where $\vec a_1'$ is the minimal element of $[\vec a_1]$ and $\vec a_2'$ is the minimal element of $[\vec a_2]$. This order is definable from $<$, so we can translate any order-invariant $\FO$ sentence on $I(\mfK)$ back to an order-invariant sentence on $\mfK$. For this backwards translation, the equality of interpreted elements is translated by $\varepsilon$, relation atoms are translated by $\rho_R$, and quantifiers in the translation are relativized to $D_I(A)$.

    Now, if local-order-invariant $\FO$ is not closed under interpretations, then certainly it does not have the same expressive power as order-invariant $\FO$ by the above. Assume instead that it is closed under interpretations. Take any $\mfK\subseteq \mfF_\sigma$. Define the $2$-dimensional interpretation $H$ between $\mfK$ and $\mfF_{\sigma\cup \{E,P\}}$ with $P$ unary and $E$ binary in which
    \begin{enumerate}
        \item $\delta(x,y)$ is a tautology,
        \item $\varepsilon(\vec x,\vec y)$ says
        \begin{align*}
            (x_1=x_2\land y_1=y_2\land x_1=y_1)\lor (x_1\neq x_2\land y_1\neq y_2),
        \end{align*}
        \item the relation $R^{H(A)}([\vec a_1],\dots,[\vec a_{a_R}])$ holds if and only if each $\vec a_i$ is of the form $(a,a)$ and $A\models R((a_1)_1,\dots,(a_{a_R})_1)$.
        \item $P^{H(A)}([\vec a])$ holds if and only if $a_1\neq a_2$, and
        \item $E^{H(A)}([\vec a],[\vec a'])$ holds if and only if $(a_1=a_2\land a'_1\neq a_2')$ holds or $(a_1\neq a_2\land a'_1=a_2')$ holds.
    \end{enumerate}
    Effectively, we have added a hub node $h$ to $A$ which is colored $P$ and is adjacent via the symmetric graph relation $E$ to every element of $A$. Via this $E$-adjacency, any order-invariant $\FO$ sentence on $\mfK$ can be translated to a local-order-invariant $\FO$ sentence on $H(\mfK)$ which uses the linear order on the substructure of $H(A)$ isomorphic to $A$ given by the local order on $S_1^{H(A)}(h)$. Because we have assumed closure under interpretations, then this local-order-invariant $\FO$ sentence on $H(\mfK)$ can be translated to a local-order-invariant $\FO$ sentence on $\mfK$.
\end{proof}

The above techniques also allow us to give a sufficient condition for the collapse of order-invariant $\FO$ to $\FO$ on classes of bounded tree-width and bounded path-width.

\begin{definition}
    A tree decomposition of a graph $G$ is a tree $T$ and a collection $X_{t_1},\dots,X_{t_n}$ indexed by the vertices of $T$ of sets of vertices of $G$ such that
    \begin{enumerate}
        \item the union $\bigcup_i X_{t_i}$ is all of $G$;
        \item for each $v\in G$, the set $S_v\subseteq T$ of vertices $t_i$ such that $v\in X_{t_i}$ forms a connected subtree of $T$; and
        \item for every pair $v,w\in G$ which are adjacent in $G$, there is some $X_{t_i}$ containing both $v$ and $w$.
    \end{enumerate}
    The width of a tree decomposition is one less than the size of the largest set $X_{t_i}$, and the tree-width of a graph $G$ is the minimum possible width of a tree decomposition of $G$. We say that a class $\mfK\subseteq \mfF_\sigma$ is of bounded tree-width if there is some $k$ such that $\mcG(A)$ is of tree-width at most $k$ for every $A\in \mfK$. Replacing the tree $T$ with a path $P$ in the above definition yields the definition of a path decomposition of $G$, and this naturally gives definitions of path-width and classes of bounded path-width.
\end{definition}

\begin{proof}[Proof of \Cref{cor:collapseonboundedtreewidth}]
    If $A$ has tree-width at most $d$, then $H(A)$ has tree-width at most $d+1$, with $H$ our hub interpretation. An order-invariant $\FO$ sentence on a class $\mfK\subseteq \mfF_\sigma$ of tree-width at most $d$ can be translated to a local-order-invariant $\FO$ sentence on the class $H(\mfK)$ of tree-width at most $d+1$, which is equivalent on $H(\mfK)$ to some $\FO$ sentence via the hypothesis of the corollary, which can be pulled back to $\mfK$ via $H$. For path-width, the same argument applies.
\end{proof}

\begin{proof}[Proof of \Cref{thm:boundedshrubdepthcollapse}]
    Let the interpretation $I$ from \Cref{thm:OPPT} consist of formulae $\delta(x)$ and $\rho_E(x,y)$. We can translate $\varphi$ to a sentence $\varphi'$ in the vocabulary $\sigma\cup\{<\}$ such that $(B,<)\models \varphi'$ if and only if $(I(B),<\harpoon D_I(B))\models \varphi$. The order-invariance of $\varphi$ on finite graphs transfers to the order-invariance of $\varphi'$ on $\mfF_\sigma$. By \Cref{thm:EEH}, there is some $\varphi''$ in the vocabulary $\sigma$ which is equivalent to $\varphi'$ on $\mfF_\sigma^{\mathrm{td},d}$. Let $U_1,\dots,U_n$ be the parameters of the transduction $T$ from \Cref{thm:OPPT}, and let $T(G,\vec U)$ be the image of $(G_k,S,\vec U)$ under the interpretation $I_T$ from the transduction $T$. We can define an $\MSO$ formula $\gamma(\vec U)$ which says that
    \begin{enumerate}
        \item $T(G,\vec U)$ has tree-depth at most $d$,
        \item the interpreted domain $\delta(T(G,\vec U))$ is exactly the elements marked by $P$ corresponding to $G$, and
        \item for all $x,y$ in $G$, the interpreted edge formula $\rho_E(x,y)$ holds in $T(G,\vec U)$ exactly when $E(x,y)$ holds in $G$. 
    \end{enumerate}
    The definitions of \textbf{(2)} and \textbf{(3)} are straightforward. For \textbf{(1)}, we define the $\MSO$ formulae $\varphi_{\mathrm{td},d}(X)$ inductively for each $d$, which will say that the induced subgraph with domain $X$ has tree-depth at most $d$. First, $\varphi_{\mathrm{td},0}(X)$ will simply say that $X$ is empty. Next, there is a $\MSO$ formula $\varphi_C(X,Y)$ which says that $Y$ is a connected component of the induced subgraph $X$, an $\MSO$ formula $\eta(Y,Y',x)$ which says that as sets $Y'=Y\setminus\{x\}$, and we define 
    \begin{align*}
        \varphi_{\mathrm{td},d+1}(X)\equiv \forall Y(\varphi_C(X,Y)\rightarrow\exists x(Y(x)\land \exists Y'(\eta(Y,Y',x)\land \varphi_{\mathrm{td},d}(Y')))).
    \end{align*}
    Translate $\varphi''$ through backwards through $T$ to an $\MSO$ sentence $\varphi_T''(\vec U)$, and define
    \begin{align*}
        \phi\equiv \exists U_1,\dots,U_n(\gamma(\vec U)\land \varphi''_T(\vec U)).
    \end{align*}
    If $G\models \varphi$, we can choose $\vec U$ such that $T(G,\vec U)\models \varphi''$. On the other hand, given $\vec U$ such that $\gamma(\vec U)\land \varphi_T''(\vec U)$, then $G\models \varphi$. Thus $\phi$ is equivalent on $\mfK$ to $\varphi$, and $\phi$ is equivalent on $\mfK$ to some first-order $\psi$ by \Cref{thm:ChenFlum}.
\end{proof}

\begin{proof}[Proof of \Cref{rem:nouniformencoding}]
    We say that a formula $\varphi(\vec x,\vec y)$ is stable on a class $\mfK\subseteq \mfK_\sigma$ if there is an integer $k\in\omega$ such that for every $A\in\mfK$, there do not exist tuples $\vec a_1,\dots,\vec a_k$ of elements of $A$ of the same length as $\vec x$ and $\vec a_1',\dots,\vec a_k'$ of elements of $A$ of the same length as $\vec y$ such that for all $i,j\in \{1,\dots,k\}$, we have $A\models \varphi(\vec a_i,\vec a_j')$ if and only if $i<j$. We say that $\mfK$ is stable if every formula is stable on $\mfK$. Roughly, this says that we cannot interpret linear orders of arbitrary length in elements of $\mfK$. 
    
    Clearly, then, the class of finite linear orders is not stable. On the other hand, every graph class of bounded shrub-depth is stable (this is a consequence of the results in \cite{ganianEtAl2012} that was first observed in \cite{gajarskyHlineny2015}). Stability is preserved under first-order interpretations, in the sense that if $\mfK$ is stable then $I(\mfK)$ is stable. The Gaifman graph of every finite linear order is complete, and so the class of finite linear orders has bounded shrub-depth. If an encoding as described in the remark existed, the class of the colored graphs encoding linear orders would be of bounded shrub-depth. However, this violates the preservation of stability under interpretations.
\end{proof}

\end{document}

%% file: preamble.tex
\newcommand{\mfF}{\mathfrak{F}}
\newcommand{\mcG}{\mathcal{G}}
\newcommand{\mfJ}{\mathfrak{J}}
\newcommand{\mfK}{\mathfrak{K}}

\newcommand{\Inc}{\operatorname{Inc}}
\newcommand{\FO}{\mathrm{FO}}
\newcommand{\MSO}{\mathrm{MSO}}
\newcommand{\CFO}{\mathrm{CFO}}
\newcommand{\tp}{\operatorname{tp}}
\newcommand{\harpoon}{\upharpoonright}